\documentclass{article}%
\usepackage{graphicx}
\usepackage{amsmath}
\usepackage{amsfonts}
\usepackage{amssymb}%
\providecommand{\U}[1]{\protect\rule{.1in}{.1in}}
\newtheorem{theorem}{Theorem}[section]

\newtheorem{lemma}[theorem]{Lemma}

\newtheorem{remark}[theorem]{Remark}

{\catcode`\@=11\global\let\AddToReset=\@addtoreset
\AddToReset{equation}{section}

\AddToReset{theorem}{section}

\begin{document}

\title{Correlation between two quasilinear elliptic problems with a source term
involving the function or its gradient}
\author{Haydar Abdel Hamid\thanks{Laboratoire de Math\'{e}matiques et Physique
Th\'{e}orique, CNRS\ UMR 6083, Facult\'{e} des Sciences, 37200 Tours, France.
\textit{Email address}:abdelham@lmpt.univ-tours.fr}
\and Marie Fran\c{c}oise Bidaut-V\'{e}ron\thanks{Laboratoire de Math\'{e}matiques
et Physique Th\'{e}orique, CNRS\ UMR 6083, Facult\'{e} des Sciences, 37200
Tours, France. \textit{Email address}: veronmf@univ-tours.fr}}
\date{.}
\maketitle

\textbf{Abstract}

Thanks to a change of unknown we compare two elliptic quasilinear problems
with Dirichlet data in a bounded domain of $\mathbb{R}^{N}.$ The first one, of
the form $-\Delta_{p}u=\beta(u)\left\vert \nabla u\right\vert ^{p}+\lambda
f(x),$ where $\beta$ is nonnegative, involves a gradient term with natural
growth. The second one, of the form $-\Delta_{p}v=\lambda f(x)(1+g(v))^{p-1}$
where $g$ is nondecreasing, presents a source term of order $0$. The
correlation gives new results of existence, nonexistence and multiplicity for
the two problems.\medskip

\textbf{R\'{e}sum\'{e}}

\textbf{Corr\'{e}lation entre deux probl\`{e}mes quasilin\'{e}aires
elliptiques avec terme de source relatif \`{a} la fonction ou \`{a} son
gradient }A l'aide d'un changement\textbf{ }d'inconnue\textbf{ }nous comparons
deux probl\`{e}mes elliptiques quasilin\'{e}aires avec conditions de Dirichlet
dans un domaine born\'{e} $\Omega$ de $\mathbb{R}^{N}.$ Le premier, de la
forme $-\Delta_{p}u=\beta(u)\left\vert \nabla u\right\vert ^{p}+\lambda f(x),$
o\`{u} $\beta$ est positif, comporte un terme de gradient \`{a} croissance
critique. Le second, de la forme $-\Delta_{p}v=\lambda f(x)(1+g(v))^{p-1}$
o\`{u} $g$ est croissante, contient un terme de source d'ordre 0. La
comparaison donne des r\'{e}sultats nouveaux d'existence, nonexistence et
multiplicit\'{e} pour les deux probl\`{e}mes.\medskip\ 

{\Large Version fran\c{c}aise abr\'{e}g\'{e}e}\medskip

Soit $\Omega$ un domaine born\'{e} r\'{e}gulier de $\mathbb{R}^{N}(N\geqq2)$
et $1<p<N$. Dans cette note nous comparons deux probl\`{e}mes
quasilin\'{e}aires. Le premier comporte un terme de source d'ordre $1$:%
\begin{equation}
-\Delta_{p}u=\beta(u)\left\vert \nabla u\right\vert ^{p}+\lambda
f(x)\hspace{0.5cm}\text{dans }\Omega,\qquad u=0\qquad\text{sur }\partial
\Omega, \tag{1}%
\end{equation}
o\`{u} $\beta$ $\in C^{0}(\left[  0,L\right)  ),$ $L\leqq\infty,$ \`{a}
valeurs $\geqq0,$ $\lambda>0$ et $f\in L^{1}(\Omega),$ $f\geqq0$ p.p. dans
$\Omega.$ Le second probl\`{e}me comporte un terme de source d'ordre $0$:
\begin{equation}
-\Delta_{p}v=\lambda f(x)(1+g(v))^{p-1}\hspace{0.5cm}\text{dans }\Omega,\qquad
v=0\qquad\text{sur }\partial\Omega, \tag{2}%
\end{equation}
o\`{u} $g\in C^{1}(\left[  0,\Lambda\right)  ),$ $\Lambda\leqq\infty,$
$g(0)=0$ et $g$ est croissante.

Le changement d'inconnue
\[
v(x)=\Psi(u(x))=%
%TCIMACRO{\dint \nolimits_{0}^{u(x)}}%
%BeginExpansion
{\displaystyle\int\nolimits_{0}^{u(x)}}
%EndExpansion
e^{\gamma(\theta)/(p-1)}d\theta,\quad\text{o\`{u} }\gamma(t)=%
%TCIMACRO{\dint \nolimits_{0}^{t}}%
%BeginExpansion
{\displaystyle\int\nolimits_{0}^{t}}
%EndExpansion
\beta(\theta)d\theta,
\]
conduit formellement du probl\`{e}me (1) au probl\`{e}me (2), et $\beta$ et
$g$ sont li\'{e}s par la relation $\beta(u)=(p-1)g^{\prime}(v).$ En
particulier $\beta$ est croissant si et seulement si $g$ est convexe. Le
changement d'inconnue inverse formel, apparemment moins utilis\'{e}, est
donn\'{e} explicitement par
\[
u(x)=H(v(x))=%
%TCIMACRO{\dint \nolimits_{0}^{v(x)}}%
%BeginExpansion
{\displaystyle\int\nolimits_{0}^{v(x)}}
%EndExpansion
\frac{ds}{1+g(s)}.
\]
Toutefois dans la transformation peuvent s'introduire des mesures. Notons
$M_{b}^{+}(\Omega)$ l'espace des mesures de Radon positives born\'{e}es sur
$\Omega$, et $M_{s}^{+}(\Omega)$ le sous-ensemble des mesures concentr\'{e}es
sur un ensemble de $p$-capacit\'{e} $0$. Nous \'{e}tablissons une
correspondance pr\'{e}cise entre les deux probl\`{e}mes:\medskip

\textbf{Th\'{e}or\`{e}me 1 }\textit{Soit }$u$\textit{ une solution
renormalis\'{e}e du probl\`{e}me }%
\begin{equation}
-\Delta_{p}u=\beta(u)\left\vert \nabla u\right\vert ^{p}+\lambda
f(x)+\alpha_{s}\hspace{0.5cm}\text{dans }\Omega,\qquad u=0\qquad\text{sur
}\partial\Omega, \tag{3}%
\end{equation}
\textit{o\`{u} }$\alpha_{s}$\textit{ }$\in M_{s}^{+}(\Omega),$\textit{ et
}$0\leqq u(x)<L$\textit{ p.p. dans }$\Omega$\textit{. Alors il existe }%
$\mu_{s}\in M_{b}^{+}(\Omega),$\textit{ telle que }$v=\Psi(u)$\textit{ est
solution atteignable du probl\`{e}me}%
\begin{equation}
-\Delta_{p}v=\lambda f(x)(1+g(v))^{p-1}+\mu_{s}\hspace{0.5cm}\text{dans
}\Omega,\qquad v=0\qquad\text{sur }\partial\Omega. \tag{4}%
\end{equation}
\textit{R\'{e}ciproquement soit }$v$\textit{ une solution renormalis\'{e}e de
(4), telle }$0\leqq v(x)<\Lambda$\textit{ p.p. dans }$\Omega,$\textit{ o\`{u}
}$\mu_{s}\in M_{s}^{+}(\Omega)$\textit{. Alors il existe }$\alpha_{s}\in
M_{s}^{+}(\Omega)$\textit{ telle que }$u=H(v)$\textit{ est solution
renormalis\'{e}e de (3). De plus, si }$\mu_{s}=0,$\textit{ alors }$\alpha
_{s}=0.$\textit{ Si }$L=\infty$\textit{ et }$\beta\in L^{1}((0,\infty
)),$\textit{ alors }$\mu_{s}=e^{\gamma(\infty)}\alpha_{s}.$\textit{ Si
}$L<\infty,$\textit{ ou }$L=\infty$\textit{ et }$\beta\not \in L^{1}%
((0,\infty)),$\textit{ et }$\alpha_{s}\not =0$\textit{ alors (3) n'a pas de
solution. Si }$\Lambda<\infty$\textit{ et }$\mu_{s}\not =0,$\textit{ alors (4)
n'a pas de solution.}\medskip

Dans le cas $\beta$ constant, les r\'{e}sultats suivants g\'{e}n\'{e}ralisent
ceux de \cite{AAP} relatifs au cas $p=2$:\medskip

\textbf{Th\'{e}or\`{e}me 2 }\textit{On suppose que }$\beta(u)\equiv
p-1,$\textit{ donc }$v=\Psi(u)=e^{u}-1$\textit{ et }$g(v)=v,$\textit{ et que}%
\[
\lambda_{1}(f)=\inf\left\{  (\int_{\Omega}\left\vert \nabla w\right\vert
^{p}dx)/(\int_{\Omega}f\left\vert w\right\vert ^{p}dx):w\in W_{0}^{1,p}%
(\Omega)\backslash\left\{  0\right\}  \right\}  >0
\]
\textit{ Si }$\lambda>\lambda_{1}(f),$\textit{ ou }$\lambda=\lambda_{1}%
(f)$\textit{ et }$f\in L^{N/p}(\Omega),$\textit{ alors (1) et (2) n'ont pas de
solution renormalis\'{e}e.}

\noindent\textit{Si }$0<\lambda<\lambda_{1}(f)$\textit{ alors (2) a une
solution unique }$v_{0}\in W_{0}^{1,p}(\Omega)$\textit{, et (1) a une solution
unique }$u_{0}\in W_{0}^{1,p}(\Omega)$\textit{ telle que }$e^{u_{0}}-1\in
W_{0}^{1,p}(\Omega)$\textit{. Si de plus }$f\in L^{r}(\Omega)$\textit{ avec
}$r>N/p,$\textit{ alors }$u_{0}$\textit{ et }$v_{0}\in L^{\infty}(\Omega
);$\textit{ et pour toute mesure }$\mu_{s}\in M_{s}^{+}(\Omega)$\textit{, (4)
a une solution renormalis\'{e}e }$v_{s}$\textit{, et donc (1) a une
infinit\'{e} de solutions }$u_{s}=H(v_{s})\in W_{0}^{1,p}(\Omega)$\textit{
moins r\'{e}guli\`{e}res que }$u_{0}$\textit{.}\medskip

Le th\'{e}or\`{e}me 1 et l'utilisation du probl\`{e}me (1) nous permettent de
d\'{e}duire un r\'{e}sultat important pour le probl\`{e}me (2), \'{e}tendant
un r\'{e}sultat classique de \cite{BrCMR} dans le cas $p=2:$\medskip

\textbf{Th\'{e}or\`{e}me 3 }\textit{On suppose que }$\Lambda=\infty,$ \textit{
}$\lim_{s\rightarrow\infty}g(s)/s=\infty,$\textit{ }$g$\textit{ est convexe
\`{a} l'infini, et }$f\in L^{r}(\Omega)$\textit{ avec }$r>N/p$\textit{. Alors
il existe }$\lambda^{\ast}>0$\textit{ tel que pour tout }$\lambda\in\left(
0,\lambda^{\ast}\right)  $\textit{ le probl\`{e}me (2) a une solution minimale
born\'{e}e }$\underline{v}_{\lambda}$\textit{, et pour tout }$\lambda
>\lambda^{\ast}$\textit{ il n'a aucune solution renormalis\'{e}e. }\medskip

Nous \'{e}tudions aussi les propri\'{e}t\'{e}s de la fonction extr\'{e}male
$v^{\ast}=\sup_{\lambda\nearrow\lambda^{\ast}}\underline{v}_{\lambda}$
\'{e}tendant certains r\'{e}sultats de \cite{CaSa}, \cite{Ne}, \cite{Sa}. Dans
le cas o\`{u} $g$ est \`{a} croissance limit\'{e}e par l'exposant de Sobolev,
nous obtenons des r\'{e}sultats d'existence d'une seconde solution
variationnelle, nouveaux m\^{e}me dans le cas $p=2,$ \'{e}tendant ceux de
\cite{AAP} et de \cite{Fe}\textbf{.}

\section{Introduction and main results}

Let $\Omega$ be a smooth bounded domain in $\mathbb{R}^{N}(N\geqq2)$ and
$1<p<N$. In this Note we compare two quasilinear problems. The first one
presents a source gradient term with a natural growth:%
\begin{equation}
-\Delta_{p}u=\beta(u)\left\vert \nabla u\right\vert ^{p}+\lambda
f(x)\hspace{0.5cm}\text{in }\Omega,\qquad u=0\hspace{0.5cm}\text{on }%
\partial\Omega, \label{PU}%
\end{equation}
where $\beta$ $\in C^{0}(\left[  0,L\right)  ),$ $L\leqq\infty,$ and $\beta$
is nonnegative, and $\lambda>0$ is a given real, and $f\in L^{1}(\Omega),$
$f\geqq0$ a.e. in $\Omega.$ Here $\beta$ can have an asymptote at point $L,$
and is not supposed to be nondecreasing.

The second problem involves a source term of order $0,$ with the same
$\lambda$ and $f:$%
\begin{equation}
-\Delta_{p}v=\lambda f(x)(1+g(v))^{p-1}\hspace{0.5cm}\text{in }\Omega,\qquad
v=0\hspace{0.5cm}\text{on }\partial\Omega, \label{PV}%
\end{equation}
where $g\in C^{1}(\left[  0,\Lambda\right)  ),$ $\Lambda\leqq\infty,$ $g(0)=0$
and $g$ is nondecreasing.

Problems of type (\ref{PU}) and (\ref{PV}) have been intensively studied the
last twenty years. The main questions are existence, according to the
regularity of $f$ and the value of $\lambda,$ regularity and multiplicity of
the solutions, and existence with possible measure data.

It is well known that the change of unknown in (\ref{PU})
\[
v(x)=\Psi(u(x))=%
%TCIMACRO{\dint \nolimits_{0}^{u(x)}}%
%BeginExpansion
{\displaystyle\int\nolimits_{0}^{u(x)}}
%EndExpansion
e^{\gamma(\theta)/(p-1)}d\theta,\text{ where }\gamma(t)=%
%TCIMACRO{\dint \nolimits_{0}^{t}}%
%BeginExpansion
{\displaystyle\int\nolimits_{0}^{t}}
%EndExpansion
\beta(\theta)d\theta,
\]
leads formally to problem (\ref{PV}), where $\Lambda=\Psi(L)$ and $g$ is given
by $g(v)=e^{\gamma(u)/(p-1)}-1.$ This is a way for studying problem (\ref{PU})
from the knowledge of problem (\ref{PV}). It is apparently less used the
reverse correlation, even in case $p=2$: for \textbf{any} function $g$
nondecreasing on $\left[  0,\Lambda\right)  $, the substitution in (\ref{PV})
\[
u(x)=H(v(x))=%
%TCIMACRO{\dint \nolimits_{0}^{v(x)}}%
%BeginExpansion
{\displaystyle\int\nolimits_{0}^{v(x)}}
%EndExpansion
\frac{ds}{1+g(s)}%
\]
leads formally to problem (\ref{PU}), where $\beta$ is defined on $\left[
0,L\right)  $ with $L=H(\Lambda)$; indeed $H=\Psi^{-1}.$ And $\beta$ is linked
to $g$ by relation $\beta(u)=(p-1)g^{\prime}(v).$ In particular $\beta$ is
nondecreasing if and only if $g$ is convex; and $L$ is finite if and only if
$1/(1+g)\not \in L^{1}\left(  0,\Lambda\right)  .$ \medskip

\noindent\textbf{Some examples} \textbf{with} $p=2$.

1. $\beta(u)=1$ and $1+g(v)=1+v.$

2. $\beta(u)=q/(1+(1-q)u),$ $q\in\left(  0,1\right)  ,$ and $1+g(v)=(1+v)^{q}%
.$

3. $\beta(u)=1+e^{u}$ and $1+g(v)=(1+v)(1+\ln(1+v)).$

4. $\beta(u)=q/(1-(q-1)u),$ $q>1$ and $1+g(v)=(1+v)^{q}.$

5. $\beta(u)=1/(1-u)$ and $1+g(v)=e^{v}.$

6. $\beta(u)=q/(1-(q+1)u),\;$ $q>0$ and $1+g(v)=1/(1-v)^{q}.\smallskip
$\textit{\ }

It had been observed in \cite{FeMu2} that the correspondence between $u$ and
$v$ is more complex, because some measures can appear, in particular in the
equation in $v$. Our first main result is to make precise the link between the
two problems. We denote by $M_{b}(\Omega)$ the set of bounded Radon measures,
and by $M_{s}(\Omega)$ the subset of measures concentrated on a set of
$p$-capacity $0$. And $M_{b}^{+}(\Omega)$ and $M_{s}^{+}(\Omega)$ are the
positive cones of $M_{b}(\Omega),$ $M_{s}(\Omega)$, and $M_{0}(\Omega)$ is the
subset of measures absolutely continuous with respect to the $p$-capacity.
Recall that $M_{b}(\Omega)=M_{0}(\Omega)+M_{s}(\Omega).$

We recall one definition of renormalized solutions among four of them given in
\cite{DMOP}. Let $U$ is measurable and finite a.e. in $\Omega$, such that
$T_{k}(U)$ belongs to $W_{0}^{1,p}(\Omega)$ for any $k>0.$ One still denotes
by $u$ the cap$_{p}$-quasi-continuous representative. Let $\mu=\mu_{0}+\mu
_{s}^{+}-\mu_{s}^{-}\in M_{b}(\Omega)$. Then $U$ is a renormalized solution of
problem
\begin{equation}
-\Delta_{p}U=\mu\hspace{0.5cm}\text{in }\Omega,\qquad U=0\quad\text{on
}\partial\Omega,\label{mu}%
\end{equation}
if $\left\vert \nabla U\right\vert ^{p-1}{\in}L^{\tau}(\Omega),$ {for any
}$\tau\in\left[  1,N/(N-1)\right)  ,$ and for any $k>0,$ there exist
$\alpha_{k},\beta_{k}\in M_{0}(\Omega)\cap M_{b}^{+}(\Omega),$ concentrated on
the sets $\left\{  U=k\right\}  $ and $\left\{  U=-k\right\}  $ respectively,
converging in the narrow topology to $\mu_{s}^{+},\mu_{s}^{-}$ such that for
any $\psi\in W_{0}^{1,p}(\Omega)\cap L^{\infty}(\Omega),$
\[
\int_{\Omega}\left\vert \nabla T_{k}(U)\right\vert ^{p-2}\nabla T_{k}%
(U).\nabla\psi dx=\int_{\left\{  \left\vert U\right\vert <k\right\}  }\psi
d\mu_{0}+\int_{\Omega}\psi d\alpha_{k}-\int_{\Omega}\psi d\beta_{k}.
\]

\begin{theorem}
\label{TP} Let $u$ be any renormalized solution of problem
\begin{equation}
-\Delta_{p}u=\beta(u)\left\vert \nabla u\right\vert ^{p}+\lambda
f(x)+\alpha_{s}\hspace{0.5cm}\text{in }\Omega,\qquad u=0\hspace{0.5cm}\text{on
}\partial\Omega, \label{PA}%
\end{equation}
where $\alpha_{s}$ $\in M_{s}^{+}(\Omega)$ and such that $0\leqq u(x)<L$ a.e.
in $\Omega$. Then there exists $\mu_{s}\in M_{b}^{+}(\Omega),$ such that
$v=\Psi(u)$ is a reachable solution of problem%
\begin{equation}
-\Delta_{p}v=\lambda f(x)(1+g(v))^{p-1}+\mu_{s}\hspace{0.5cm}\text{in }%
\Omega,\qquad v=0\hspace{0.5cm}\text{on }\partial\Omega. \label{PS}%
\end{equation}

Conversely let $v$ be any renormalized solution of (\ref{PS}), where $\mu
_{s}\in M_{s}^{+}(\Omega)$ and such that $0\leqq v(x)<\Lambda$ a.e. in
$\Omega$. Then there exists $\alpha_{s}\in M_{s}^{+}(\Omega)$ such that
$u=H(v)$ is a renormalized solution of (\ref{PA}).

Moreover if $\mu_{s}=0,$ then $\alpha_{s}=0.$ If $L=\infty$ and $\beta\in
L^{1}((0,\infty)),$ then $\mu_{s}=e^{\gamma(\infty)}\alpha_{s}.$ If $L<\infty$
or if $L=\infty$ and $\beta\not \in L^{1}((0,\infty)),$ and $\alpha_{s}%
\not =0,$ then (\ref{PA}) has no solution. If $\Lambda<\infty$, and $\mu
_{s}\not =0,$ then (\ref{PS}) has no solution.
\end{theorem}

This theorem extends in particular the results of \cite{AAP} where $p=2$ and
$L=\infty$. The nonexistence result when $\beta\not \in L^{1}((0,\infty)),$
and $\alpha_{s}\not =0,$ answers to an open question of \cite{Por}.

First we apply to the case $\beta$ constant, which means $g$ linear.

\begin{theorem}
\label{T2}Assume that $\beta(u)\equiv p-1,$ thus $v=\Psi(u)=e^{u}-1$ and
$g(v)=v.$ Suppose that
\begin{equation}
\lambda_{1}(f)=\inf\left\{  (\int_{\Omega}\left\vert \nabla w\right\vert
^{p}dx)/(\int_{\Omega}f\left\vert w\right\vert ^{p}dx):w\in W_{0}^{1,p}%
(\Omega)\backslash\left\{  0\right\}  \right\}  >0 \label{lam}%
\end{equation}
If $\lambda>\lambda_{1}(f),$ or $\lambda=\lambda_{1}(f)$ and $f\in
L^{N/p}(\Omega),$ then (\ref{PU}) and (\ref{PV}) admit no renormalized solution.

\noindent If $0<\lambda<\lambda_{1}(f)$ there exists a unique solution
$v_{0}\in W_{0}^{1,p}(\Omega)$ to (\ref{PV}), thus a \textbf{unique} solution
$u_{0}\in W_{0}^{1,p}(\Omega)$ to (\ref{PU}) such that $e^{u_{0}}-1\in
W_{0}^{1,p}(\Omega)$. If $f\in L^{r}(\Omega)$ with $r>N/p,$ then $u_{0}$ and
$v_{0}\in L^{\infty}(\Omega),$ and moreover for any measure $\mu_{s}\in
M_{s}^{+}(\Omega)$, there exists a renormalized solution $v_{s}$ of
(\ref{PS}); then there exists an \textbf{infinity} of less regular solutions
$u_{s}=H(v_{s})\in W_{0}^{1,p}(\Omega)$ of (\ref{PU}).
\end{theorem}

\begin{remark}
Under the assumption (\ref{lam}), most of these existence results extend to
general $g$ such that $\Lambda=\infty$ and $\lim\sup_{\tau\longrightarrow
\infty}g(\tau)/\tau<\infty$. They extend to the case
\begin{equation}
\lim\sup_{\tau\longrightarrow\infty}g(\tau)/\tau^{q}<\infty\quad\text{for some
}q\in\left(  1,N/(N-p)\right)  \label{maja}%
\end{equation}
if moreover $f\in L^{r}(\Omega)$ with $qr^{\prime}<N/(N-p).$
\end{remark}

Next consider problem (\ref{PV}) with general $g,$ and $f\in L^{r}(\Omega)$
with $r>N/p.$ It is easy to prove that for small $\lambda>0$ there exists a
minimal solution $\underline{v}_{\lambda}\in W_{0}^{1,p}(\Omega)$ such that
$\left\Vert \underline{v}_{\lambda}\right\Vert _{L^{\infty}(\Omega)}<\Lambda$.
Our main result is an extension of a well known result of \cite{BrCMR} for
$p=2,$ and \cite{CaSa} for $p>1$. It is noteworthy that \textbf{the proof uses
problem }(\ref{PU}):

\begin{theorem}
\label{truc}Assume that $\Lambda=\infty,$ and $\lim_{s\longrightarrow\infty
}g(s)/s=\infty$, and $g$ is convex near infinity, and $f\in L^{r}(\Omega)$
with $r>N/p$. There exists a real $\lambda^{\ast}>0$ such that:

(i) for $\lambda\in\left(  0,\lambda^{\ast}\right)  $ problem (\ref{PV}) has a
minimal \textbf{bounded} solution $\underline{v}_{\lambda}$,

(ii) for $\lambda>\lambda^{\ast}$ there exists no \textbf{renormalized} solution.
\end{theorem}

When $g$ is subcritical with respect to the Sobolev exponent $p^{\ast
}=Np/(N-p)$, we obtain \textbf{new multiplicity results} for problem
(\ref{PV}), even in the case $p=2,$ extending \cite{AAP} and \cite{Fe} :

\begin{theorem}
\label{main} Under the assumptions of Theorem \ref{truc}, assume that
\begin{equation}
\lim\sup_{\tau\longrightarrow\infty}g^{p-1}(\tau)/\tau^{Q}<\infty\text{ for
some }Q\in\left(  1,p^{\ast}-1\right)  , \label{majet}%
\end{equation}
and $f\in L^{r}(\Omega)$ with $(Q+1)r^{\prime}<p^{\ast}.$ Then there exists
$\lambda_{0}>0$ such that for any $\lambda<\lambda_{0},$ there exists at least
\textbf{two solutions} $v\in W_{0}^{1,p}(\Omega)\cap L^{\infty}(\Omega)$ of
(\ref{PV}). Moreover if $p=2,$ $g$ is convex, or $g$ satisfies the
Ambrosetti-Rabinowitz growth condition and $f\in L^{\infty}(\Omega)$, the same
result holds with $\lambda_{0}=\lambda^{\ast}.$
\end{theorem}

Concerning the extremal solution, we get the following, extending some results
of \cite{CaSa}, \cite{Sa}:

\begin{theorem}
\label{limi}Under the assumptions of Theorem \ref{truc}, the extremal function
$v^{\ast}=\sup_{\lambda\nearrow\lambda^{\ast}}\underline{v}_{\lambda}$ is a
renormalized solution of (\ref{PV}) with $\lambda=\lambda^{\ast}$. If
$N<p(1+p^{\prime})/(1+p^{\prime}/r)$ then $v^{\ast}$ $\in W_{0}^{1,p}%
(\Omega).$ Moreover $v^{\ast}$ $\in L^{\infty}(\Omega)$ in any of the
following conditions:

(i) $N$ is arbitrary and (\ref{majet}) holds and $(Q+1)r^{\prime}<p^{\ast},$

(ii) $N$ is arbitrary and (\ref{maja}) holds and $qr^{\prime}<N/(N-p),$

(iii) $N<pp^{\prime}/(1+1/(p-1)r).$
\end{theorem}

\begin{remark}
Using Theorems \ref{TP}, \ref{truc} and \ref{main}, we deduce existence and
nonexistence results for problem \ref{PU}. In Theorem \ref{TP}, function $f$
can depend on $u$ or $v,$ which strongly extends the range of applications.
For example, taking $g(v)=v,$ and $f=u^{b},$ $b>0,$ problem $-\Delta
_{p}u=(p-1)\left\vert \nabla u\right\vert ^{p}+\lambda u^{b}$ relative to $u$
leads to $-\Delta_{p}v=\lambda(1+v)^{p-1}\ln^{b}(1+v)$ relative to $v.$ Then
for small $\lambda$ the problem in $u$ has an infinity of solutions $u\in
W_{0}^{1,p}(\Omega),$ two of them being bounded.
\end{remark}

\begin{remark}
A part of our results is based on a growth assumption on $g$. Returning to
problem (\ref{PU}), this condition is not always easy to verify. When
$L=\infty$, all the "usual" functions $\beta,$ even with a strong growth,
satisfy $\lim\sup_{\tau\longrightarrow\infty}g(\tau)/\tau^{q}<\infty$ for any
$q>0,$ see \cite{AAP}. However using the converse correlation between $g$ and
$\beta,$ we prove that the conjecture that this condition always holds is
\textbf{wrong}: let $F\in C^{0}(\left[  0,\infty\right)  )$ be any strictly
convex function, with $\lim_{s\longrightarrow\infty}F(s)=\infty.$ Then there
exists an increasing function $\beta$ such that $\lim_{t\longrightarrow\infty
}\beta(t)=\infty$ and the corresponding $g$ satisfies $\lim\sup_{\tau
\longrightarrow\infty}g(\tau)/F(\tau)=\infty.$
\end{remark}

\section{Sketch of the main proofs}

In some proofs we use a regularity Lemma:

\begin{lemma}
\label{boot}Let $1<p<N,$ and $F\in L^{m}(\Omega),$ and $\overline
{m}=Np/(Np-N+p)\ ($thus $1<\overline{m}<N/p).$ Let $U$ be a renormalized
solution of problem
\[
-\Delta_{p}U=F\quad\text{in }\Omega,\qquad U=0\quad\text{on }\partial\Omega.
\]
If $1<m<N/p,$ then $U^{(p-1)}\in L^{k}(\Omega),$ with $k=Nm/(N-pm).$ If
$m=N/p,$ then $U^{(p-1)}\in L^{k}(\Omega)$ for any $k\geqq1.$ If $m>N/p,$ then
$U\in L^{\infty}(\Omega).$ If $1<m<\overline{m},$ then $\left\vert \nabla
U\right\vert ^{(p-1)}\in L^{\tau}(\Omega),$ with $\tau=Nm/(N-m).$ If
$m\geqq\overline{m},$ then $U\in W_{0}^{1,p}(\Omega).$
\end{lemma}

\textbf{Proof of Theorem \ref{TP} }For $p\neq2,$ we cannot use approximations
of the equations because of the nonuniqueness of the solutions of $-\Delta
_{p}U=\mu$ with $\mu\in M_{b}^{+}(\Omega)$. The main idea is to use the
equations satisfied in the sense of distributions by the truncated functions
$T_{K}(u)=\min(K,u)$ and $T_{k}(v)=\min(k,v)$ with $k=\Psi(K)$, using
definition (ii) of renormalized solution given above: \
\begin{align*}
-\Delta_{p}T_{K}(u)  &  =\beta(T_{K}(u))\left\vert \nabla\,T_{K}(u)\right\vert
^{p}+\lambda f\chi_{\left\{  u\leqq K\right\}  }+\alpha_{K},\text{ \qquad in
}\mathcal{D}^{\prime}\left(  \Omega\right)  ,\\
-\Delta_{p}T_{k}(v)  &  =\lambda f(1+g(v))^{p-1}\chi_{\left\{  v\leqq
k\right\}  }+\mu_{k},\text{\qquad in }\mathcal{D}^{\prime}\left(
\Omega\right)  ,
\end{align*}
where $\mu_{k}$ and $\alpha_{K}$ are two measures concentrated on the same
set: $\left\{  u=K\right\}  =\left\{  v=k\right\}  ,$ and explicitely related
by $\mu_{k}=(1+g(k))^{p-1}\alpha_{K},$ and respectively converging weakly * to
$\mu_{s}$ and $\alpha_{s}$. The nonexistence results are consequences of some
properties of renormalized solutions, also called Inverse Maximum
Principle.\smallskip

\textbf{Proof of Theorem \ref{T2}. }The nonexistence is first proved for
(\ref{PU}), and then for (\ref{PV}) by Theorem \ref{TP}. The existence is
obtained by iteration and approximation, using \cite{DMOP}. Uniqueness follows
from Picone's identity, adapted to renormalized solutions.\smallskip

\textbf{Proof of Theorem \ref{truc}. }Formally, if $v$ is a solution of
(\ref{PV}) for some $\lambda,$ and $u=H(v),$ then $\bar{u}=(1-\varepsilon)u$
is a supersolution of (\ref{PU}) relative to $\bar{\lambda}=(1-\varepsilon
)^{p-1}\lambda,$ and $\bar{v}=\Psi(\bar{u})$ is a supersolution of (\ref{PV})
relative to $\bar{\lambda};$ then there exists a solution $v_{1}\leqq\bar{v}.$
Using Theorem \ref{TP} we show that it is not formal, since actually
\textbf{no measure appears}. In the (best) case $H\left(  \infty\right)
<\infty,$ $\bar{v}$ is bounded, then also $v_{1}$ is bounded. Otherwise a
bootstrapp using Lemma \ref{boot}is needed for constructing a bounded
solution.\smallskip

\textbf{Proof of Theorem \ref{main}.} The Euler function $J_{\lambda}$ is well
defined, and for small $\lambda$ it has the geometry of Mountain Path near 0,
but the Palais-Smale sequences may be unbounded$.$ From \cite{Je}, (\ref{PV})
has a second solution for almost any small $\lambda,$ and then for a sequence
$\lambda_{n}\rightarrow\lambda,$ and we are lead to prove that the solutions
$v_{\lambda_{n}}$ relative to $\lambda_{n}$ converge to a solution to
(\ref{PV}) relative to $\lambda$. The usual estimates for the case $p=2,$
using an eigenfunction as test function, cannot be extended. We get an
estimate of $-\Delta_{p}v_{\lambda_{n}}$ in $L^{1}(\Omega)$ in another way,
using the convexity of $g$. The estimate of $v_{\lambda_{n}}$ in $W_{0}%
^{1,p}(\Omega)$ is obtained by contradiction. For larger $\lambda,$ if $p=2$,
$J_{\lambda}$ has still the geometry of mountain path near $\underline
{v}_{\lambda}$; the question is open when $p\neq2$. Under
Ambrosetti-Rabinowitz condition we apply some results of \cite{GhPr}%
.\smallskip

\textbf{Proof of Theorem \ref{limi}. }The estimates come from Lemma \ref{boot}
and well known regularity results for quasilinear equations, and from an
extension of techniques of \cite{Ne}.

\end{document}